\documentclass{article}
\usepackage[top=0.825in,bottom=0.675in,left=0.6in,right=0.4in, paperwidth=8.5in, paperheight=11in]{geometry}
\usepackage{CBI}
\usepackage{newtxmath}
\title{Simplifying and Refactoring Introductory Calculus}
\author{Jonathan Bartlett}

\begin{document}
\maketitle

\begin{abstract}
First year calculus is often taught in a way that is very burdensome to the student.
Students have to memorize a diversity of processes for essentially performing the same task.
However, many calculus processes can be simplified and streamlined so that fewer concepts can provide more flexibility and capability for first-year students.
\end{abstract}

\section{Introduction}

While the exact set of topics in any particular calculus book or course may vary, the general method of calculus training has been essentially set in stone for the last hundred years.
Nonetheless, there are many issues with this methodology that have been insufficiently addressed over the years.
These issues range from the sequencing of topics to the content of the topics themselves.

A well-recognized phenomena in computer programming is known as ``code debt.''
Code debt occurs when new ideas and changes get incorporated into a program, but the rest of the program doesn't change sufficiently to take these new ideas and changes into account.
Because of this, there winds up being a lot of duplication and confusion over the right ways of doing things.
It is known as code \emph{debt} because eventually, to relieve the tension, future work will have to be done to the rest of the system to bring it back into alignment.

To relieve code debt, computer programmers engage in in a process known as ``refactoring.''
The idea behind refactoring is to reconceptualize the whole of the computer program in order to discover which facilities are the core, distinguishable pieces and which ones are merely a variant or permutation of those pieces.

The present goal of the paper is to begin a similar process of refactoring the subject of calculus and the way that it is taught.

\section{Proof Mechanisms and Learning Calculus}

Calculus books and courses are often written so that each step of the development of calculus can be proved.
Most books begin with limits, and then calculus is then proved using limits.
Unfortunately, most calculus students are not yet ready for limits, and limits often wind up being a confusing side-challenge that derail students before they have even begun to study the subject in depth.

In nearly any subject, students do best when proofs occur \emph{after} the core content is learned.
As an example, consider children learning a language.
Children learn English (or any other language) long before they learn the grammar of English.
In fact, most students can speak English perfectly well without ever learning the rules (grammar) that govern English.
Even when learning foreign languages, immersive approaches (i.e., those in which the language precedes the grammar of it) tend to work most effectively \citep{tlbt2015}.

Of course, mathematical learning is not identical to language learning.
Part of mathematics \emph{is} proof, so learning proof tools is certainly an important part of mathematics.
Additionally, with most technical subjects, having a grounding in the ``why'' aspects certainly helps learn the intuitions behind the ``how.''
However, even when the proofs help build an intuition, there is a difference between teaching an intuition behind an idea by proving its truth and trying to get students to construct their own proofs ex nihilo.

Therefore, the sequencing in this method will focus on first making the intuitive understanding work, and then only later establish proofs.

The primary target of this reasoning is on limits.
Most of the time, limits are introduced \emph{before} any of the problems it is intended to solve, which makes the subject confusing in the first place.
If a student has grown up doing the simple calculation $f(3)$ when $f(x) = x^2$, it makes no sense to transform this into limit notation, $\lim\limits_{x\to 3} x^2$.
This is simply seen by the student as extra steps which have no perceptible benefit.

Students see the new notation, a bunch of rules, and are asked to put them together to solve trivial problems that they could solve using the normal rules of algebra.
Therefore, they start to get lost because they are bombarded with extra notation for nothing.
They can't grasp at what this could possibly be for, and are lost before they ever really start.

At minimum, the concept of limits should be \emph{introduced} to solve the problem of having ``holes'' in a graph (rather than for expressions that are easily calculated without limits).
This at least provides a good context for why we should bother with them.
In my own classes, I have moved limits to the end of the course altogether.
Some want to introduce limits to establish a formal definition of continuity.
However, it is much more straightforward to new students to introduce continuity as ``being able to draw a graph without picking up your pencil'' than to talk about limits.
Students understand drawing without picking up your pencil much better than limits. 
After a whole year of talking about continuity in terms of whether or not you pick up your pencil (and why that is important to calculus), you can say that, as another benefit of limits, we get a more formal definition of what it means to have a continuous graph.

Thus, the explicit justification for the derivative occurs long after the students learn the derivative and learn to use it well.
This has the additional advantage of including L'Hospital's rule as part of the general discussion of limit behavior.

\section{Introducing the Derivative}
\label{secderivintro}

In this method, derivatives are introduced through extensive practical work on evaluating slopes between two points on a graph.
Given a graph of, say, $y = x^2$, what is the slope between the points $x = 1$ and $x = 2$?
Basic algebra is used to to find the slope from two points, using the standard formula $m = \frac{y_1 - y_0}{x_1 - x_0}$.

Next, we discuss what it means to find the slope \emph{at} a point.
We talk about the problems of finding a slope at a single point (the formula we have reduces to $\frac{0}{0}$), and then discuss methods of getting around that.
I suggest estimating the slope at the point using two points that are close together.
Several problems are worked for getting the distance between the two $x$ locations closer and closer ($x = 1$ and $x = 1.1$, then move to $x = 1$ and $x = 1.0001$).

Eventually, we try to make a formula that expresses the ``general idea'' of a slope between two points on a given line.
If we are looking at two points where the $x$ values are $0.01$ away from each other, we eventually come up with the formula $m = \frac{(x + 0.01)^2 - x^2}{(x + 0.01) - x}$ to match their existing knowledge of the slope formula.

After solving a few of these, I then suggest that we introduce a parameter to represent the distance.
This way, we can decide \emph{later} how close we want the points to be.
If we call our parameter $h$, then the formula becomes $m = \frac{(x + h)^2 - x^2}{(x + h) - x}$.
This simplifies to $m = 2x + h$.
So, to find the slope between any two points that are $h = 0.01$ away from each other, we simply use the formula.
We then make $h$ smaller and smaller and smaller to get a ``more exact'' value of the slope at a point.
I then ask them to think about if there is any way to abuse this formula in order to find a formula for the slope at a single point.
At this point, students are ready to understand that if I set $h$ to $0$ that the slope is equivalent to the slope at a particular point.

I then generalize this with functions, using $f(x)$ to represent the function we want the slope of, and show that basic simplifications give us the formulat for the derivative, $y' = \frac{f(x + h) - f(x)}{h}$.
I then say that, for any particular $f()$, after getting rid of $h$ in the denominator, we can simply substitute $0$ for $h$.  

I do mention that $h$ is not really zero, but in fact really close to it.
Nonetheless, I tell them that the formal justification for treating $h$ as zero will come later in the course, and that for now they should treat it as a value that is sufficiently close to zero that it can be treated as such, but sufficiently far away from zero that it doesn't produce divide-by-zero errors.
Doing this helps students start building intuitions around these types of situations, which will help us refine them in our formal discussions of limits later on.

\section{Derivatives vs. Differentials}

There are two potential operators to use for teaching calculus.
The first, and more commonly used, is the derivative operator, which can be written as $D_x()$ or as $\ddx()$.
This is commonly spoken of as ``taking the derivative with respect to $x$.''

The problem with this operator is that it becomes very confusing for students when they get to derivatives of implicit functions.
The reason for this is that, when applied to implicit functions, derivatives are spitting out other derivatives rather than simple answers.
And, in fact, it does this for some variables and not others.
For instance, the derivative $x$ with respect to $x$ is $1$, but the derivative of $y$ with respect to $x$ is $\frac{\dy}{\dx}$.  This assymetry causes confusion.

While the derivative operation is perfectly consistent and solves the necessary problems, a much clearer way to teach students all types of derivatives is to instead focus on the differential operator, $\diffop()$.
The differential operator is similar to the derivative, but it does not treat any variable specially.
There is no variable that the derivative is being taken with respect to.
Instead, whichever variable(s) are in the expression always comes out at the end as a differential.

To understand the differences, we will apply the operators to both an explicit and an implicit function.
The explicit function will be $y = x^3$ and the implicit function will be $xy = 5$.

Using the $\ddx()$ operator, the explicit function is straightforward.
\begin{align*}
\ddx(y) &= \ddx(x^3) \\
\dydx &= 3x^2.
\end{align*}
However, the implicit function is unintuitive for students.
\begin{align*}
\ddx(xy) &= 5 \\
x\,\dydx + y &= 0
\end{align*}
The product rule, which is supposed to have some amount of symmetry to it, now seems to have no symmetry whatsoever.

However, when taking differentials, the operations are much more straightforward.
First, the explicit function.
\begin{align*}
\diffop(y) &= \diffop(x^3) \\
\dy &= 3x^2\,\dx
\end{align*}
As you can see, this is identical to the $\ddx()$ operation except that the result is not divided by $\dx$.
This localizes the differentials to the contexts where they occur.
As we will see shortly, this simplifies a lot of reasoning about complicated derivatives.

Here is the differentiation of the implicit function.
\begin{align*}
\diffop(xy) &= \diffop(5) \\
x\,\dy + y\,\dx &= 0
\end{align*}
This is functionally identical to the implicit derivative.
However, it maintains the symmetry of the product rule itself.
Additionally, it required no extra rules for the student for the handling of $y$ vs.\ $x$.
They are both treated \emph{identically} in the differentiation step.

To convert the result of differentiation into a derivative, the student merely needs to solve for $\dydx$ using normal algebra rules.
Therefore, \emph{all} derivatives become a two-step process:
\begin{enumerate}
\item Find the differential of the equation.
\item Solve for the derivative you are interested in.
\end{enumerate}

This way, the process for finding $\dydx$ is the same as for finding $\dxdy$.
There are no funny steps.
You merely differentiate both sides and then solve for the particular derivative (ratio of differentials) that you want.
It does not matter if the equation is implicit, explicit, or, as we will see shortly, multivariate, the process is identical.

Even though going through the differential is slightly more complicated for simple equations, overall I have found that students find the separation of finding differentials and then finding derivatives to be much easier to process over a wider variety of problems.

Interestingly, this is the way the original Leibnizian calculus was developed.
The process focused on differentials rather than derivatives, and did not focus on finding ratios of differentials or worrying about identifying independent variables \citep{bos1974}.
As mentioned in Section~\ref{secderivintro}, I do \emph{start} with the traditional derivative, as it is much more concrete.
However, I do attempt to pivot very quickly from derivatives to differentials.

Another possible mechanism for doing this was pioneered by \citet{thompson1910}.
In this method, differentials are introduced first, basically as almost-discrete differences.
Essentially, the question is asked, given, say, $y = x^2$, if I added a number named $\dy$ to all the $y$s, I would have to compensate by adding some other number $\dx$ to all the $x$s.
This possibility is also explored, but not very fully, by \citet{dray2010}.
Future work may entail discovering which approach is more beneficial to students---starting with slopes as a concrete tie-in to algebra and then pivoting to differentials, or starting more directly with differentials at the very beginning.

\section{Treating Differentials Algebraically}

Because we are solving for the derivative instead of calculating it directly, we must be willing to treat differentials algebraically.
Although in the 19th century this was frowned upon, developments in the 20th century have shown that differentials can be thought of as algebraic units.
Abraham Robinson's formal description of the hyperreals and the introduction of non-standard analysis provides sufficient justification for treating differentials algebraically \citep{robinson1974,henle2003,keisler2012}.
Essentially, differentials are treated as hyperreal values---values which exist on the number line but which are not contained by the reals.

Although differentials can be treated algebraically, they can't actually be solved for as specific numbers.
Eventually, to get a number that can be evaluated or calculated, the differentials have to be in ratio with each other.
However, for intermediate steps, there is no problem with treating differentials such as $\dx$ and $\dy$ exactly as you would the variables $x$ and $y$.
They can be multiplied, divided, cancelled, etc.
Treating differentials as algebraic units is very intuitive for students who have studied and practiced the algebraic treatment of unknowns for quite some time.

Additionally, the ability to treat differentials algebraically improves the comprehension of topics such as related rate problems.
In related rate problems, $\dt$ springs up as if from nowhere.
Books often teach this as ``taking the derivative with respect to $t$,'' but this is confusing since there is no $t$.
I have found that students are much less confused by an algebraic approach.
We simply take the differential the same way we always have, but then we divide both sides of the equation by $\dt$.
This works simply because, algebraically, we can divide by any value we wish as long as we divide by the same thing on both sides of the equation.

Again, by using differentials instead of derivatives, we have transformed a number of processes that students find unintuitive into a single process where the intuition is supplied by the student's knowledge of algebra.

\section{Multivariable Calculus}
\label{SecMultivariable}

This method of using differentials instead of derivatives also simplifies some aspects of multivariable calculus.
Take the equation $z^2 = xy$.
This can be easily converted into differentials.
\begin{align*}
\diffop(z^2) &= \diffop(xy) \\
2z\,\dz &= x\,\dy + y\,\dx 
\end{align*}
This can then be solved algebraically for any derivative that is desired.
For instance, $\frac{\dy}{\dz}$ can be found.
\begin{align*}
2z\,\dz &= x\,\dy + y\,\dx  \\
2z\,\dz - y\,\dx &= x\,\dy \\
\frac{2z}{x}\dz - \frac{y}{x}\dx &= \dy \\
\frac{2z}{x} - \frac{y}{x}\frac{\dx}{\dz} &= \frac{\dy}{\dz}
\end{align*}
What this means is that the rate at which $y$ changes with respect to $z$ not only depends on the actual values of $x$, $y$, and $z$, but that it \emph{also} depends on the rate at which $x$ changes with respect to $z$.
Many calculus books skip over this idea, probably because when differentials are not treated independently, these processes are very complicated.
However, with dealing with differentials instead of derivatives, these ideas arise very naturally from the notation.

\section{Partial Derivatives/Differentials}

This treatment of differentials also leads to a very straightforward way of defining partial differentials and derivatives.
Total derivatives define the relationships between changes in all parts of the system.
But what happens if you want to hold some pieces still, and only find the way that one variable influences another when treated by itself?

In order to do that, we would say that the other variables \emph{don't change}.
Another way of saying that these variables don't change is to say that their differentials are zero.

In other words, to convert from a multivariable differential (as in Section~\ref{SecMultivariable}) to a partial differential is to merely set the other differentials which aren't under consideration to zero.

For instance, Section~\ref{SecMultivariable} ended with the equation
\begin{align*}
\frac{\dy}{\dz} &= \frac{2z}{x} - \frac{y}{x}\frac{\dx}{\dz}.
\end{align*}
To understand the relationship between the changes in $y$ and $z$ \emph{if everything else is kept constant}, we merely need to set the differentials of everything else ($\dx$ in this case) to zero.
Doing this results in the equation
\begin{align*}
\frac{\dy}{\dz} &= \frac{2z}{x} - \frac{y}{x}\frac{\dx}{\dz} \\
\frac{\partial y}{\partial z} &= \frac{2z}{x} - \frac{y}{x}\frac{0}{\dz} \\
 &= \frac{2z}{x} - 0 \\
 &= \frac{2z}{x}
\end{align*}
Therefore, the partial derivative of $y$ with respect to $z$ is $\frac{2z}{x}$.

Note, however, that partial differentials cannot always be treated algebraically. 
This is not due to a failing of the concept of partial derivatives, but merely of their notation.
For our original equation, $\partial y$ can refer to two \emph{different} entities in the ratios $\frac{\partial y}{\partial z}$ (where $\dx = 0$) and $\frac{\partial y}{\partial x}$ (where $\dy = 0$).
Information about \emph{which} particular $\partial y$ is being spoken about is contained in the denominator, and therefore splitting the numerator from the denominator results in a loss of information.\footnote{A possible solution to this problem would be to subscript partial differentials with the list of differentials which were allowed to change.  
  This is cumbersome, but allows for an algebraic treatment of partial differentials.  
  So, for instance, $\frac{\partial y}{\partial z}$ would be written as $\frac{\partial_{yz}y}{\partial_{yz}z}$ and $\frac{\partial y}{\partial x}$ would be written as $\frac{\partial_{yx}y}{\partial_{yx}x}$.  
  Doing this makes it clear $\partial_{yz}y$ is a distinct algebraic entity from $\partial_{yx}y$. 
  This is the subject of a current paper in progress from the present author.}

\section{Higher Order Differentials and Derivatives}

Performing derivatives by taking differentials first also has benefits down the road, as it allows students to understand other notations more fully.
The Liebniz notation for the second derivative $\left(\frac{\diffop^2 y}{\dx^2}\right)$ has long baffled many students.
Most books take a ``just use it and don't ask questions'' approach.\footnote{\citet{bartlett2018} gives several examples of textbooks taking this approach.}
However, recent advances have shown that the standard notation for higher-order derivatives is not only baffling, it is in a very real sense incorrect \citep{bartlett2018}.

To understand the issues, first recognize that differentials are actually a shorthand.
When you take the differential of a composite function, you always wind up with a differential of the inner function.
For instance, 
$$ \diffop(\sin(x^2)) = \cos(x^2)\,\diffop(x^2) = \cos(x^2)\,2x\,\dx. $$
However, the last term, $\dx$, stands in the same relationship to $2x$ as $\diffop(x^2)$ does to $\cos(x^2)$. 
It is the differential of the interior function.
Therefore, $\dx$ actually is just a shorthand for $\diffop(x)$.
Since it is irreducible, it is shortened to just $\dx$. 
This shortening happens both to reduce reading and writing effort, and also to mark the fact that this differential cannot be further reduced.
Nonetheless, $\dx$ is actually a composite---an operator and an operand.\footnote{To emphasize this, I usually typeset $\dx$ such that the $\diffop$ is in roman type and the $x$ is in italic.  This is similar to the way that other functions such as $\sin x$ are typeset.}

But what do $\diffop^2y$ and $\dx^2$ mean?
The latter is straightforward enough.  
$\dx^2$ is simply a shorthand for $(\diffop(x))^2$.
However, $\diffop^2y$ (with the superscript after the differential operator) actually means applying the differential operator twice.
In other words, while $\dy$ is short for $\diffop(y)$, $\diffop^2y$ is a shorthand for $\diffop(\diffop(y))$.

To see how this plays out, imagine the equation 
\begin{equation}
\label{EqHodExample}
y = x^3.
\end{equation}
Let's start by taking the differential of this equation twice.
\begin{align}
y &= x^3 \nonumber \\
\diffop(y) &= \diffop(x^3) \nonumber \\
\dy &= 3x^2\,\dx && \text{first differential} \nonumber \\
\nonumber \\
\diffop(\dy) &= \diffop(3x^2\,\dx) \nonumber \\
\hdiff{y}{2} &= 3x^2\,\hdiff{x}{2} + 6x\,\dx^2 && \text{second differential} \label{EqHodSecDiff}
\end{align}
This result may seem surprising, but allow for an explanation.
The term $3x^2\,\hdiff{x}{2}$ seems like it is out of place, but it is not.
Since $3x^2\,\dx$ is the \emph{product} of $3x^2$ and $\dx$, the product rule has to be used to resolve the next differential.

Therefore, since $\diffop(uv) = u\,\dv + v\,\du$, while one of the outputs is the anticipated $6x\,\dx^2$, there is also another one, $3x^2\,\hdiff{x}{2}$.
If $x$ is the independent variable, this term goes to zero because $\hdiff{x}{2}$ goes to zero (see \citet{bartlett2018} and \citet{bos1974} for an explanation of why).
However, if $x$ is not an independent variable, the term is vitally important.
Without keeping the term, the second differential would not be algebraically manipulable.
Keeping the term then allows it to remain manipulable.

Now, let us divide the whole thing by $\dx^2$.
Doing so yields
\begin{equation}
\label{notquitesecondderiv}
\frac{\hdiff{y}{2}}{\dx^2} = 3x^2\,\frac{\hdiff{x}{2}}{\dx^2} + 6x 
\end{equation}

However, this has more terms than what we normally expect from a second derivative.
This is because the second derivative comes from the following sequence of steps:
\begin{enumerate}
\item Take a differential
\item Divide by $\dx$
\item Take another differential
\item Divide by $\dx$
\end{enumerate}
Because the derivative is the combination of steps 1 \& 2, the second derivative simply repeats those steps.
However, if you take the differential twice in a row before dividing by $\dx^2$, then you wind up with a strange-looking answer.

A few things to note about this:
\begin{enumerate}
\item The left-hand side has the form that we normally associate with the second derivative.
\item The right-hand side does not have the form that we normally associate with the second derivative.  It has an extra term in it.
\item The reason for this is that, even though it is normally associated with the second derivative, the left-hand side \emph{is not the correct notation for the derivative of the derivative of $y$}.
\item In other words, the equation given by Equation~\ref{notquitesecondderiv} is correct as far as it goes, but it is not the second derivative, because the notation we have come to associate with the second derivative actually refers to a different quantity altogether.
\end{enumerate}

According to \citet{bartlett2018}, the full notation for the second derivatives should be
\begin{equation}
\label{EqSecDeriv}
\frac{\hdiff{y}{2}}{\dx^2} - \frac{\dy}{\dx}\frac{\hdiff{x}{2}}{\dx^2}.
\end{equation}
This can be deduced simply from taking two derivatives of $y$.
The first derivative is obviously $\frac{\dy}{\dx}$.
The second derivative is found by taking the differential of $\frac{\dy}{\dx}$ and then dividing by $\dx$.
\begin{align}
\frac{\diff\left(\frac{\dy}{\dx}\right)}{\dx} &= \frac{\frac{\dx\,\diffop(\dy) - \dy\,\diffop(\dx)}{\dx^2}}{\dx} \nonumber \\
  &= \frac{\dx\,\diffop(\dy) - \dy\,\diffop(\dx)}{\dx^3} \nonumber \\
  &= \frac{\dx\,\diffop(\dy)}{\dx^3} - \frac{\dy\,\diffop(\dx)}{\dx^3} \nonumber \\
  &= \frac{\diffop(\dy)}{\dx^2} - \frac{\dy}{\dx}\frac{\diffop(\dx)}{\dx^2}
\end{align}
By noting that $\diffop(\dy) = \diffop(\diffop(y)) = \hdiff{y}{2}$ and $\diffop(\dx) = \diffop(\diffop(x)) = \hdiff{x}{2}$ we can see that this is equivalent with Equation~\ref{EqSecDeriv}.

So how does this square with Equation~\ref{notquitesecondderiv}?

If we subtract $3x^2\,\frac{\hdiff{x}{2}}{\dx^2}$ from both sides of Equation~\ref{notquitesecondderiv}, it will result in
\begin{align*}
\frac{\hdiff{y}{2}}{\dx^2} - 3x^2\,\frac{\hdiff{x}{2}}{\dx^2} &= 6x.
\end{align*}
If you recognize $3x^2$ as being the first derivative (i.e., $\frac{\dy}{\dx}$) it is apparent that the left-hand side of this equation is the same as the improved form of the second derivative listed in Equation~\ref{EqSecDeriv}.

This formula can be \emph{algebraically} rearranged to yield the second derivative of $x$ with respect to $y$, or combined with other derivative formulas (say, the derivative of $x$ with respect to $t$) to algebraically accomplish a change of variables.
Previously, this was only available using specialty formulas such as Fa\`a di Bruno's formula.

Thus, teaching using Leibnizian differentials confers numerous advantages for higher order differentials:
\begin{enumerate}
\item The notation is clearer, because there is a definitive reason for the way that the notation looks.
\item The notation allows differentials to be modified algebraically, which was not previously possible for higher order derivatives.
\item Because the notation is algebraically manipulable, the notation allows for students to easily find relationships that previously required memorized formulas.
\end{enumerate}

\section{Getting Rid of Logarithmic Differentiation}

Logarithmic differentiation is the process taught by most calculus textbooks for taking the derivative of functions of the form $u^v$.
Essentially, what is taught is to use logarithms to remove the exponent, and then take the derivative now that the exponent is removed.
The problem with this is that it is needlessly complicated, and forces the student to use different processes depending on the derivative in question.\footnote{The process of logarithmic differentiation is necessary in the \emph{proof} of some important rules, but not in their \emph{usage}.}

Before learning logarithmic differentiation, the student had a more-or-less unified process for taking derivatives: 
\begin{enumerate}
\item Look at the form of the function.
\item Find the corresponding rule (memorized or from the book).
\item Apply the rule.
\end{enumerate}
However, logarithmic differentiation breaks that process, adding extra steps in some circumstances.

If the mathematics required this, that would be one thing.
However, there is actually a rule available for forms of the type $u^v$ which is rarely mentioned even in the appendices of most calculus books.
The rule is
\begin{equation}
\label{EqGenPowerRule}
\diffop(u^v) = vu^{v-1}\du + \ln(u)u^v\dv
\end{equation}
This formula can be derived in many ways.
The most straightforward is to first set $z = u^v$ and then differentiate both sides using logarithmic differentiation (logarithmic differentiation is useful to \emph{prove} the formula for $u^v$ but after that it is fairly useless because you can just apply the formula).
Doing this yields
\begin{align*}
z &= u^v \\
\ln(z) &= \ln(u^v) \\
\ln(z) &= v\,\ln(u) \\
\diffop(\ln(z)) &= \diffop(v\,\ln(u)) \\
\frac{1}{z}\dz &= \frac{v}{u}\du + \ln(u)\dv \\
\dz &= \frac{zv}{u}\du + z\,\ln(u)\dv \\
\dz &= \frac{v\,u^v}{u}\du + \ln(u)\,u^v\dv \\
\dz &= v\,u^{v-1}\du + \ln(u)\,u^v\dv 
\end{align*}

This is a decent proof for first-year calculus.  
However, a more interesting proof can be found by simply recognizing that a total differential is the sum of its partials.
So, the partial differential of $u^v$ with $v$ kept constant is $v\,u^{v-1}\du$ and the partial differential of $u^v$ with $u$ kept constant is $\ln(u) u^v\dv$.
Therefore, the total differential is merely the sum of these, as you can see in Equation~\ref{EqGenPowerRule}.

Probably due to its rarity of actually appearing in textbooks, this rule has been given a variety of names, the two most common being the generalized power rule and the functional power rule.

Logarithmic differentiation does have some unique uses, but these are mostly upper-level ideas.
For instance, logarithmic differentiation can be used to convert products of sequences (which are difficult to integrate) into a sum of sequences (which tend to be easier).

For instance, logarithmic differentiation allows us to say that
\begin{align}
\text{If~~~~} &y = \prod_{j = 1}^n f(x, j) \nonumber \\
\text{Then~~~~} &y' = y \sum_{j = 1}^n \diffop(\ln(f(x, j)))
\end{align}

Logarithmic differentiation is also sometimes used to simplify complicated fractions, exponents, and even products, but, since there are already rules for all of these, the ``simplification'' usually just makes learning the process more difficult for students.
Students need unified processes more than they need tricks to make things easier.
Conceptual simplicity is usually preferable to speed of computation.

\section{Hyperreals for Limit Analysis}
\label{seclimits}

$\epsilon$-$\delta$ proofs have long been the bane of calculus students.
While they do present an interesting mathematical technique, the retention rate for understanding $\epsilon$-$\delta$ proofs is so low as to hardly be worth doing \citep{katz2017}.
Instead, the method which helps students understand the process of limits the most is to use the hyperreal number line.

On the hyperreal number line, $\epsilon$ is the unit of the infinitesimal (not the same as the $\epsilon$ in $\epsilon$-$\delta$ proofs).
To take a right-handed limit, the student merely replaces $x$ with $x + \epsilon$ everywhere it occurs in the expression.
This ``looks at'' the function immediately to the right of the point in question---an infinitely small step to the right.
The left-handed limit is found by subtracting $\epsilon$.

As an example, the expression $\frac{x^2 - 25}{x - 5}$ cannot be evaluated at $x = 5$ because it results in a zero in the denominator.  
However, any value \emph{except} $5$ will work.  
Therefore, if we bump $x$ an infinitely small amount to the right, the divide-by-zero problem will no longer exist.
Therefore, to find $\lim\limits_{x\to 5^+}\frac{x^2 - 25}{x - 5}$, we merely replace $x$ with $x + \epsilon$.
This yields
\begin{align*}
\lim\limits_{x\to 5^+}\frac{x^2 - 25}{x - 5} &= \frac{(x + \epsilon)^2 - 25}{x + \epsilon - 5} \\
 &= \frac{x^2 + 2x\epsilon + \epsilon^2 - 25}{x + \epsilon - 5} \\
 &= \frac{(5)^2 + 2(5)\epsilon + \epsilon^2 - 25}{(5) + \epsilon - 5} \\
 &= \frac{25 + 10\epsilon + \epsilon^2 - 25}{\epsilon} \\
 &= 10 + \epsilon
\end{align*}
This is infinitely close to $10$.
Therefore, the limit is $10$.

Additionally, I have found from experience that limits are best taught at the end of a year of calculus rather than at the beginning.
Teaching about hyperreal numbers and then using them in limits provides a good way to make many of the intuitions developed over the first year of calculus more concrete.

\section{The Integral as an Infinite Sum}

One more change that improves calculus for first-year students is redefining the integral from an area to an infinite sum.
This is a very subtle difference, but one that I have found to be important.
Finding the area under a curve is one particular \emph{usage} of the integral.
However, more generally, the integral is used as a tool of summation.

Using the idea of the integral as a tool of summation helps explain the usage of the integral as the area under the curve, but the converse is not true---explaining the integral as the area under the curve does not help students imagine other uses of the integral such as arc lengths and volumes of revolution.

When the integral is defined as an infinite sum of infinitely small pieces, then it is straightforward to then say that what is to the right of the integral defines what each individual small piece that we are adding together looks like.
Figure~\ref{FigSummations} lists several different sorts of things we can add together using integrals.
\begin{figure*}
\caption{Summations Using Integration}
\label{FigSummations}
\begin{center}
\begin{tabular}{c|c|c|c}
\textbf{general geometry} & \textbf{infinitesimal geometry} & \textbf{formula} & \textbf{integral} \\
\hline
 & & &  \\[-2ex]
area under the curve & ultrathin rectangles & $\text{height}\cdot\text{width}$ & $\myint y\,\dx$ \\[1ex]
arc length & line length & $\sqrt{(x_1 - x_0)^2 + (y_1 - y_0)^2}$ & $\myint \sqrt{\dx^2 + \dy^2}$ \\[1ex]
revolution about $x$ & ultrashort cylinders & $\pi\, r^2 h$ & $\myint \pi\,y^2\,\dx$ \\[1ex]
revolution about $y$ & ultrathin shells (curved boxes) & $ \text{height}\cdot\text{length}\cdot\text{thickness}$ & $\myint y\,2\pi x\,\dx$
\end{tabular}
\end{center}
\end{figure*}
In each of these, the integral represents the infinite sum of well-defined infinitely small objects.
If the integral is \emph{defined} as an infinite sum, this makes an intuitive connection for the student.
However, if the integral is \emph{defined} as the area under the curve, jumping out of this to understand how the area under the curve can be reconfigured as one of the other operations is quite confusing.

Additionally, defining the integral as an infinite sum allows for the integral to be more easily defined in multivariable situations.
The reason for this is that, as an infinite sum, the integral is adding up all of the \emph{differences} that occur into a \emph{total difference}.
This works just a straightforwardly for multivariable differentials as it does for single variable differentials.

For instance, let's say we have the equation
$$ \dz = \dy + \dx. $$
The integral of this is just
$$ z = y + x + C.$$
As an ``area under the curve'' this makes no sense.
However, as an ``infinite sum'' this makes perfect sense.
If $\dz$ represents the sum of $\dy + \dx$, then the total sum will be the difference of $y + x + C$ evaluated at two points.
So, traveling from $x = 2, y = 3$ to $x = 7, y = 4$ will give a sum total of
$$(7 + 4 + C) - (2 + 3 + C) = 11 + C - 5 - C = 6.$$

Therefore, treating the integral as an infinite sum not only helps students generalize the integral to various integration-related formulas, but it also helps students generalize the integral into a multivariable version, and similarly to complex number situations.

The area under the curve can be seen as a particular, easy-to-understand instance of infinite summation.

One other benefit of treating the integral as an infinite sum is that it makes the notation more clear.
Some introductory texts, when defining the integral as the area under the curve, simply use the integral form as a pro forma way of specifying the variable of integration.

For instance, it will often be explained like this:
\[
\myint \overbrace{x^3 + 2x^2 - 3x + 5}^{\text{formula to be integrated}} \,\underbrace{\dx}_{\text{variable of integration}}
\]
However, doing it this way needlessly restricts the usage of integration to only area under the curve, and makes the other uses of it harder to understand.

\section{Resistance to Hyperreals}

The primary reason that these methods are not more widespread is the lack of enthusiasm for the hyperreal number system.  
This comes from two sources.
The first is a historical bias against infinities and infinitesimals, with infinitesimals actually having the stronger of the negative biases.
The primary bias against infinitesimals comes from the fact that they were \emph{used} long before they were \emph{proved}.  

Before rigorous means of working with infinitesimals were established, many critics pointed to the inconsistent ways in which infinitesimals were handled as proof that they were non-entities.
In some circumstances, infinitesimals were treated as zero and thrown away, and, in other circumstances, they were used in denominators and were therefore treated as non-zero entities which could also be used for algebraic cancelling.
In the latter half of the twentieth century, Robinson's hyperreal number system made these operations rigorous with the hyperreal number system \citep{robinson1974}.
However, by that time, the damage from hundreds of years of skepticism had already taken its toll.

The most classic recrimination against infinitesimals was done by George Berkeley, in his famouse quote:
\begin{quote}
And what are these Fluxions? The Velocities of evanescent Increments? And what are these same evanescent Increments? They are neither finite Quantities nor Quantities infinitely small, nor yet nothing. May we not call them the ghosts of departed quantities? \citep{berkeley1734}\footnote{In modern terminology, \emph{fluxions} refer to derivatives, and \emph{evanescent increments} refer to infinitesimals.}
\end{quote}

Additionally, even though the hyperreal number system has been shown to be usable as a consistent system for the inclusion of infinities and infinitesimals, it is not the only viable candidate.  
Cantor's transfinite number system \citep{cantor1915}, surreal numbers \citep{knuth1974}, dual numbers \citep{wolfe2014}, and other systems have all been proposed for extending the real number line into infinitesimals.
While the hyperreals are the most widespread, it is not the only system available.
This lack of standard convention has prevented a lot of building on any one foundation.

Finally, there is also a philosophical distrust of hyperreal numbers based on general concerns over infinities and infinitesimals.
Both infinities and infinitesimals have troubled certain schools of mathematicians.
These concerns can be divided into ontological concerns and epistemological concerns.

The ontological concern is about whether infinities exist and/or are needed in mathematics.  This concern stems from the line of mathematicians following in the footsteps of Leopold Kronecker.
Kronecker opposed Cantor's theory of infinities because Kronecker only admitted mathematical concepts which could be constructed in a finite number of steps from the natural numbers \citep{dauben1990}.

David Hilbert's program was similar to Kronecker's.
However, Hilbert was not actively opposed to infinity per se, but believed that infinities were ideals and the natural numbers were reality.
Therefore, Hilbert believed that any true statement about natural numbers that were proved with infinity could also be proven without them.
Infinities weren't invalid, they were just superfluous.
Thus, Hilbert's program was to define a finite set of axioms which were consistent and complete, and could prove any valid theorem without explicitly relying on any concepts of infinity \citep{zach2016}.
G\"odel later proved that Hilbert's program was unworkable, but the motivations and concerns behind it remain today.
Even for mathematicians who accept Cantor and G\"odel's infinities, a continued concern about over-reliances on infinities and infinitesimals remains.

For the epistemological concern, I cannot point to any one bright source.  
However, this concern stands as a subtext to many conversations I have had regarding hyperreal numbers, infinities, and infinitesimals.
Essentially, the idea is that since we cannot point to anything that is infinitely small or infinitely big in the world around us, it is not safe to make solid assertions about such entities, as those assertions cannot be tested.
Our minds (and therefore our proofs) can always fail, therefore, building on premises which can only be proved logically and not physically is dangerous.

While these concerns are understandable, the long-term payoff will come from helping students see the concepts of the infinite and infinitesimal more clearly, and that will only come from practice and familiarity with the concepts.
The infinite should not be feared, in fact, mathematics is one area which allows us to grapple with the infinite on a much more rigorous basis.

\section{Conclusion}

By making a few modifications to the way that calculus is taught, students can be presented a more unified, holistic system.
This is both easier to use and easier to understand.
Having a single process and expanding out its usage to more and more complex cases is much more straightforward than having to reinvent the system at every step, and forcing students to memorize different processes for different situations.

The changes proposed here include:
\begin{itemize}
\item \textbf{Separating differentiation from finding the derivative.}  This allows the unification of explicit differentiation, implicit differentiation, and multivariable differentiation.  Additionally, it helps explain (and correct) the notions of higher-order differentials.
\item \textbf{Treating differentials algebraically.}  This is an extension of the separation of finding differentials and derivatives.  Additionally, when done correctly, it improves the usability of higher-order differentials and derivatives.
\item \textbf{Using a rule for $u^v$ instead of logarithmic differentiation.}  Instead of forcing students to use different types of processes for different forms, teaching the generalized power rule allows students to take the differential of $u^v$ directly, just as with every other form.
\item \textbf{Using hyperreals for limits.}  Hyperreals allow for a more intuitive approach to limits.  Additionally, moving limits to the end of a first-year course allows students to develop intuitions around the derivative first before seeing the formal proof of their validity.
\item \textbf{Treating integrals as infinite sums.}  Treating the integral as an infinite sum instead of an area under a curve allows for easier generalization of the concept of the integral into various geometric situations (and even non-geometric situations).  Additionally, this allows for a more straightforward generalization of integration into the inclusion of multiple variables.
\end{itemize}

The goal is to simultaneously simplify introductory calculus while making it more powerful at the same time.
This is accomplished by (using computer science terminology) ``refactoring'' calculus into pieces that are more easily recombined, adapted, and applied to various situations.
Since differentials can always be transformed into derivatives by algebraic rearrangement, nothing is lost in their treatment as individuated entities.

These ideas have been incorporated into a new text on calculus \citep{bartlett2018cgu}, and future study is needed to fully assess the impact of these ideas (both positive and negative) on the teaching of calculus.

\bibliographystyle{CBI}
\bibliography{SimplifiedCalculus}

\end{document}